\documentclass[reqno,12pt]{amsart}

\usepackage{mathrsfs, amssymb, latexsym, amsmath, graphicx, psfrag}
\newtheorem{thm}{Theorem}

\newtheorem{lemma}[thm]{Lemma}
\newtheorem{conj}{Conjecture}

\theoremstyle{remark}

\renewcommand{\le}{\leqslant}
\renewcommand{\ge}{\geqslant}
\newcommand{\SM}{\mathscr{S}}
\newcommand{\bb}{{\bf b}}
\newcommand{\bbz}{\mathbb{Z}}

\renewcommand{\a}{\alpha}
\renewcommand{\d}{\delta}
\newcommand{\tla}{T\!L^{\!A}_n}
\newcommand{\tlb}{T\!L^{\!B}_n}
\newcommand{\da}{D^A_n(\d)}
\newcommand{\db}{D^B_n(\a, \d)}
\newcommand{\A}{{\bf A}}

\newcommand{\p}[2]{\psi_{#1,#2}}
\newcommand{\g}[2]{\gamma_{{#1},{#2}}}
\newcommand{\ph}[1]{\phi_{#1}}
\renewcommand{\b}[2]{\beta_{#1,#2}}
\newcommand{\D}[2]{{\bf D}_{#1,#2}}
\newcommand{\im}{\mathrm{Im}}
\newcommand{\CL}{\mathrm{NC}}

\newcommand{\la}{\langle}
\newcommand{\ra}{\rangle}

\title{The Gram determinant of the type B Temperley-Lieb algebra}
\author{Qi Chen\ \ and \ \ Jozef H. Przytycki}

\begin{document}
\maketitle

\section{Introduction}\label{sec1}
\addtocounter{footnote}{1}

In this paper, we solve a problem posed by the late Rodica Simion
regarding type $B$ Gram determinants, cf. \cite{schmidt}.
We present this in a fashion influenced by the work of 
W.B.R.Lickorish on Witten-Reshetikhin-Turaev invariants of 3-manifolds.
We will give a history of this problem in a sequel paper in
which we also plan to address other related questions by Simion 
\cite{Sim,schmidt}
and connect the problem to Frenkel-Khovanov's work \cite{FK}.

\section{The type B Gram determinant}

\begin{figure}[ht]
\centering
\psfrag{1}{\vspace{-.3cm}\hspace{-.3cm}$a_1$}
\psfrag{2}{\hspace{-.3cm}$a_2$}
\psfrag{3}{\hspace{-.0cm}$a_n$}
\psfrag{4}{\hspace{.1cm}$a_{n+1}$}
\psfrag{5}{\hspace{-.0cm}$a_{2n}$}
\psfrag{6}{\hspace{-.0cm}$a_{2n-1}$}
\psfrag{7}{\hspace{-.0cm}$S$}
		\includegraphics[height=1.5in]{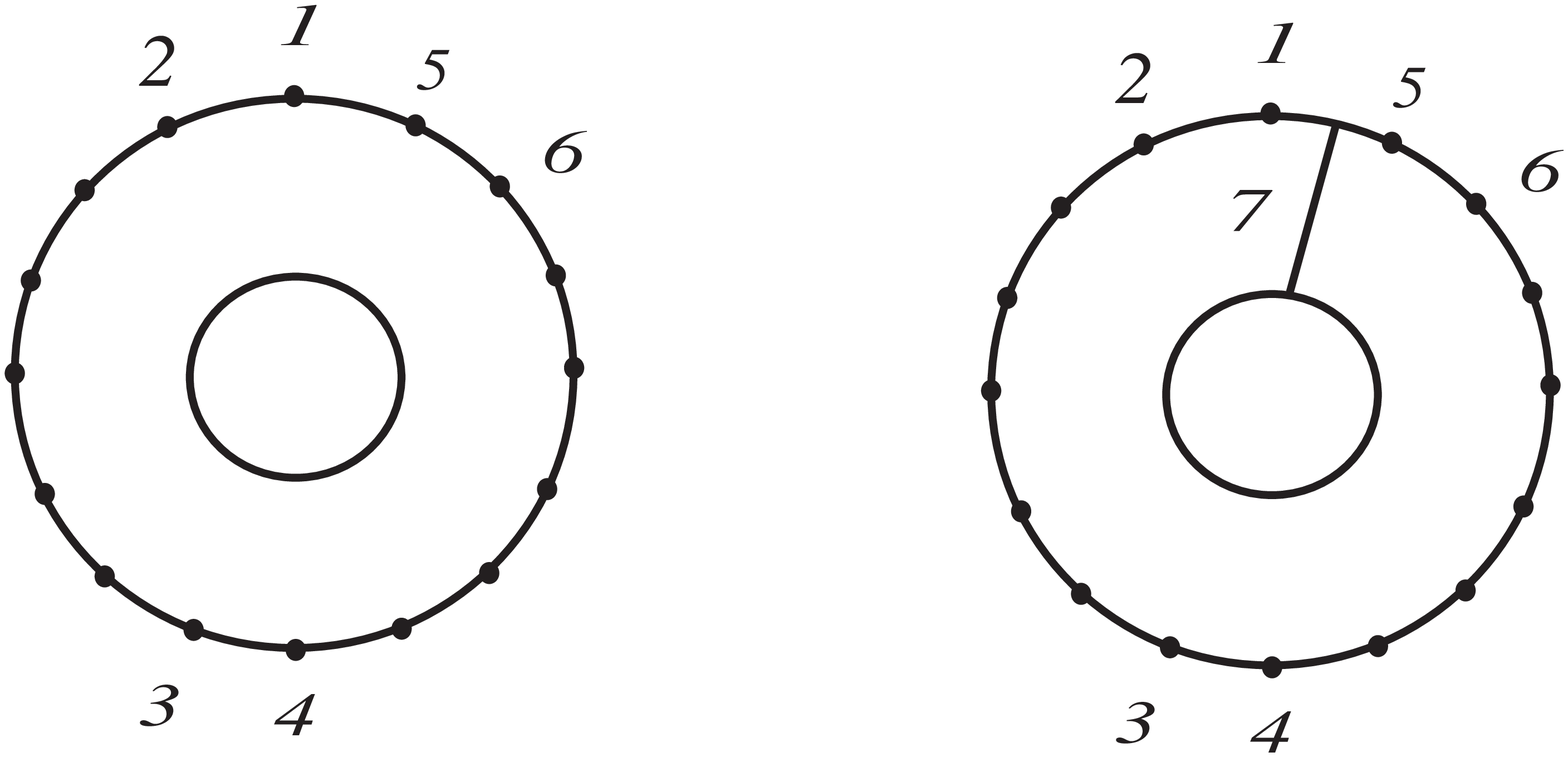} 

\caption{$\A_n$ and 
$\A_n$ with a segment $S$}
\label{fig1}
\end{figure}

\begin{figure}[ht]
	\centering
		\includegraphics[height=.7in]{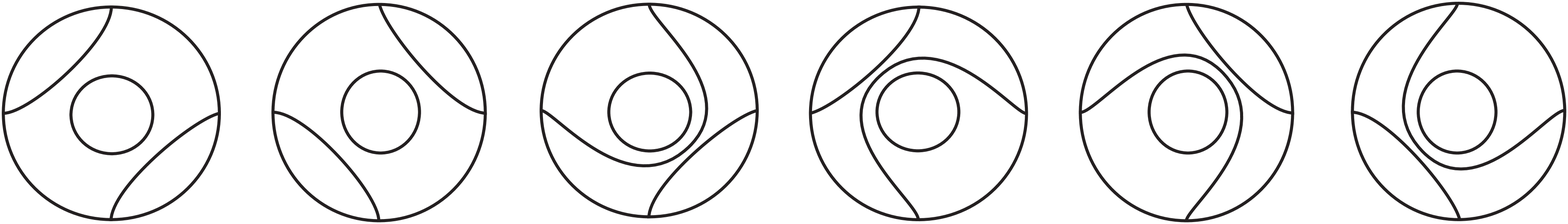}
\caption{Connections in $\bb_2 = \{b_1, b_2, b_3, b_4, b_5\}$}\label{fig2}
\end{figure}

\begin{figure}[ht]
	\centering
		\includegraphics[height=1.5in]{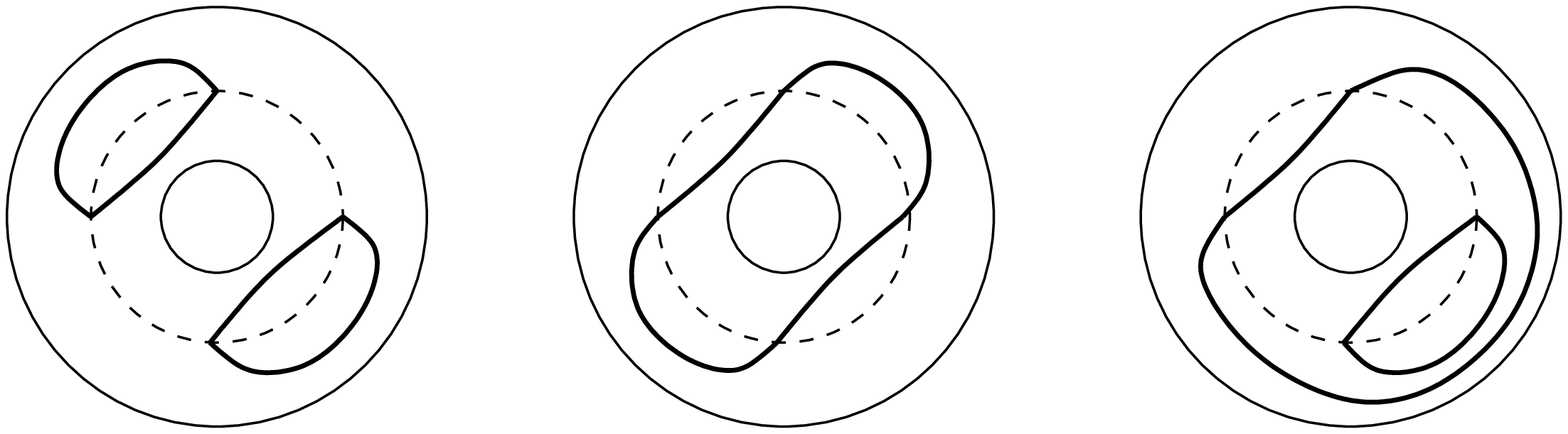}
\caption{$\la b_1, b_1\ra = \d^2$,\ $\la b_1, b_2\ra = \a$,\ 
$\la b_1, b_3\ra = \a \d$; \ \ $b_i\in \bb_2$}
\label{Ffig3}
\end{figure}

Let $\A_n$ be an annulus with $2n$ points, $a_1,\ldots, a_{2n}$, 
 on the outer circle of the boundary, cf. Fig \ref{fig1}.  
Let $\bb_n := \{b_1, b_2, \ldots, b_{\binom{2n}n}\}$ be the set of all 
possible diagrams, up to deformation, in $\A_n$ with $n$ non-crossing 
chords connecting these $2n$ points, Fig \ref{fig2}.
We define a pairing $\la\ ,\ \ra$ on $\bb_n$ as follows:
Given $b_i, b_j\in\bb_n$ we glue $b_i$ with the inversion of $b_j$ along the 
marked circle, respecting the labels of the marked points. The resulting
picture is an annulus with two types of disjoint circles, 
homotopically non-trivial and trivial; compare Fig \ref{Ffig3}. 
The bilinear form $\la\ ,\ \ra$ is define by $\la b_i, b_j\ra = \a^m \d^n$
where $m$ and $n$ are the number of homotopically non-trivial circles 
and homotopically trivial circles respectively.

Let
$$
G_n(\a,\d) = \Big( \la b_i, b_j \ra \Big)_{1\le i,j \le \binom{2n}n}
$$
be the matrix of the pairing on $\bb_n$ called the Gram matrix of 
the type $B$ Temperley-Lieb algebra. We denote
its determinant by $\db$. The roots of $\db$ were 
predicted by D{\c a}bkowski and Przytycki, and the complete 
factorization of $\db$ was conjectured by G.Barad:

\begin{conj}[G. Barad]\label{thm1}
$$
\db=\prod_{i=1}^n \left( T_i(\d)^2 - \a^2 \right)^{\binom{2n}{n-i}}
$$
where $T_i(\d)$ is the Chebyshev (Tchebycheff) polynomial of the first kind:
$$
T_0 = 2, \qquad T_1 = \d, \qquad T_i = \d\, T_{i-1} - T_{i-2}.
$$
\end{conj}

The rest of the paper is devoted to a proof of Conjecture \ref{thm1}.
It follows directly from the following two lemmas, the first of which 
is proven in Section 3.

\begin{lemma}\label{l1}
For $i\ge 1$, $\a=(-1)^{i-1}T_i(\d)$ is a zero of $\db$ of multiplicity at least $\binom{2n}{n-i}$.
\end{lemma}

\begin{lemma}\label{l2}
Let $S$ be a line segment connecting the two
boundary components of $\A_n$ such that $S$ is disjoint 
from $a_i$, $1\le i\le 2n$; see \ref{fig1}.
Let $c(b_i)$ denote the number of chords in $b_i$ that cut $S$, 
and let $P=(p_{ij})$ be
a diagonal matrix defined by $p_{ii} = (-1)^{c(b_i)}$.
Then $G_n(-\a,\d) = PG_n(\a,\d)P^{-1}$.
\end{lemma}

\begin{proof} The power of $\a$ in $\la b_i, b_j\ra$ is
 congruent to $c(b_i) + c(b_j)$ modulo $2$, thus
$$
\la b_i, b_j\ra|_{\a\mapsto -\a}= (-1)^{c(b_i) + c(b_j)}
\la b_i, b_j\ra
$$
and Lemma 2 follows.
\end{proof}


\begin{proof}[Proof of Conjecture \ref{thm1}]
According to Lemma \ref{l2}, $G_n(-\a, \d)$ and $G_n(\a,\d)$ are 
conjugate matrices. Hence
$\a=(-1)^{i}T_i(\d)$ is a zero of $\db$ of the same multiplicity as 
$\a=(-1)^{i-1}T_i(\d)$. 
Therefore, by this and Lemma \ref{l1} we have
$$
\db = p \prod_{i=1}^n \left(T_i(\d)^2 - \a^2\right)^{\binom{2n}{n-i}},
$$
for some $p\in\bbz[\a,\d]$. 
The diagonal entries in $G_n(\a,\d)$ are all equal to $\d^n$ and they are of 
highest degree in each row thus $\db$ is a monic polynomial in variable $\d$ 
of degree $n\binom{2n}{n}$.  Furthermore $T_i(\d)$ is a
monic polynomial of degree $i$. Couple these with a well known 
equality\footnote{We use ``telescoping" to get
$ 2\sum_{i=1}^n i \binom{2n}{n-i} =
2\sum _{i=1}^n n(\binom{2n-1}{n-i} - \binom{2n-1}{n-i-1}) =
2n\binom{2n-1}{n-1} = n\binom{2n}n$.}, 
$$
2\sum_{i=1}^n i \binom{2n}{n-i} = n\binom{2n}n,
$$
to conclude that $p=1$.

\end{proof}

\section{Proof of Lemma \ref{l1}}

It is enough to show that the nullity of $G_n((-1)^{i-1}T_i(\d),\d)$ 
is at least $\binom{2n}{n-i}$, which
we prove by the theory of Kauffman Bracket Skein Module (KBSM); 
see \cite{hoste-przytycki, przytycki} for the definition and properties of
KBSM. Denote the KBSM of a 3-manifold $X$ by $\SM(X)$. Let $\A$ be an annulus.
For any two elements $x,y$ in 
$\SM(\A)=\bbz[A^{\pm1},\alpha]$,\footnote{For any surface $F$, 
we write $\SM(F)$ for $\SM(F\times [0,1])$. In $\SM(\A)$, $\a$ represents 
a nontrivial curve, $1$ -- the empty curve, and $\d= -A^2 - A^{-2}$ -- 
a trivial curve. $\bb_n$ is a basis of a relative KBSM, $\SM(\A_n)$ 
as a module over $\bbz[A^{\pm1},\alpha]$.}
let $H(x,y)$ be the element 
in $\SM(S^3) = \bbz[A, A^{-1}]$ obtained by decorating the
two components of the Hopf link with $x$ and $y$. 
Denote the $k$-th Jones-Wenzl idempotent by $f_k$, 
cf. \cite{lickorish2}. Define a linear map
$$
\ph k : \SM(\A) \to \bbz[A, A^{-1}]
$$
such that $\ph k(x)=H(x,\hat{f_k})$, where $\hat{f_k}\in\SM(\A)$
is the natural closure of $f_k$. For $b_i,b_j\in \bb_n$, we will 
consider $\la b_i, b_j\ra$ as an element of 
$\SM(\A)$. If $\la b_i, b_j\ra = \a^m \d^n$ then
\begin{equation}\label{eq1}
\ph k(\la b_i, b_j\ra) = (-A^{2(k+1)}-A^{-2(k+1)})^m (-A^2-A^{-2})^n\Delta_k,
\end{equation}
where $\Delta_k = (-1)^{k}(A^{2(k+1)}-A^{-2(k+1)})/(A^2-A^{-2})$ 
is the Kauffman bracket of $\hat{f_k}$;
see page 143 of \cite{lickorish2}.
To relate the Gram matrix $G_n(\a,\d)$ to the map $\ph k$ we substitute 
$\d = -A^2-A^{-2}$ to obtain  
$T_k(\d)= (-1)^k(A^{2k}+A^{-2k})$. Let
$$
F_{n,k} = \Big(\ph {k-1}(\la b_i, b_j\ra)\Big)_{1\le i,j \le \binom{2n}n}.
$$
Then
$$
G_n((-1)^{k-1}T_k(-A^2-A^{-2}), -A^2-A^{-2}) 
= \frac{1}{\Delta_{k-1}} F_{n,k}.
$$
Therefore, Lemma \ref{l1} follows from the next lemma. 

\begin{lemma}\label{lemma4}
The nullity of $F_{n,k}$ is at least $\binom{2n}{n-k}$.
\end{lemma}

To prove Lemma \ref{lemma4} we need some
linear maps defined on $\SM(\A_n)$ and $\SM(\D nk)$, $k\ge 0$,
where $\D nk$ is the disk with $2(n+k)$ points on its boundary. 
These points are labeled counter-clockwise by
$a_1,\ldots, a_{2n}$, $l_1,\ldots, l_k$ and $u_k, \ldots, u_1$. 

\begin{figure}[ht]
\centering
\psfrag{1}{\hspace{-.2cm}$\stackrel {\p 32}{\mapsto}$}
\psfrag{2}{\hspace{-.2cm}$\stackrel {\g 32}{\mapsto}$}
\psfrag{3}{\hspace{-0cm}$f_2$}
\psfrag{4}{\hspace{-.2cm}$\stackrel {\b 32}{\mapsto}$}
		\includegraphics[height=1in]{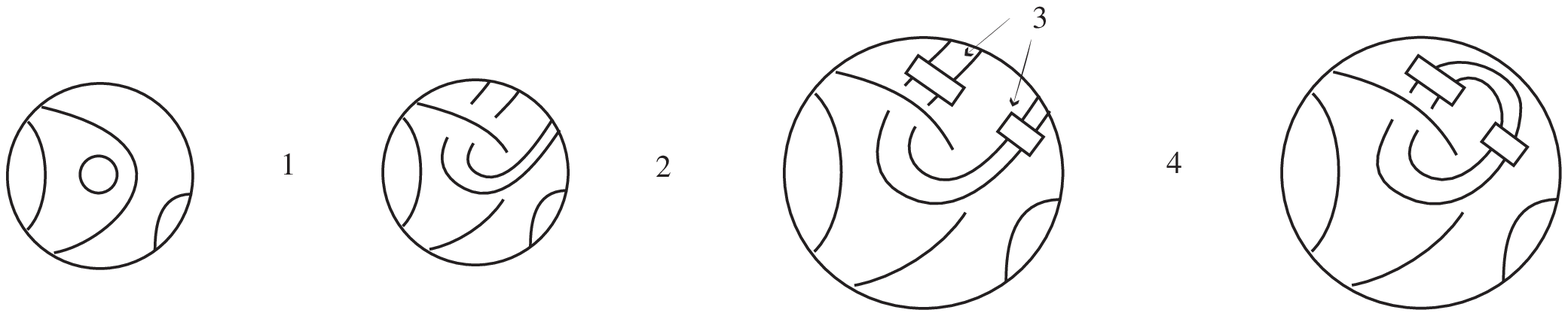}
	\caption{}\label{f7}
\end{figure}

Let
\begin{equation*}\label{eq2}
\p nk : \SM(\A_n)\to \SM(\D nk)
\end{equation*}
be a linear map defined as follows:
 Let $\p nk'$ be an embedding of $\A_n$ into a neighborhood of the 
boundary of $\D nk$ such
that the point $a_i$ on $\A_n$ is mapped to $a_i$. 
Let $L\subset\D nk$ denote the lollipop consisting of
the image of the inside boundary of $\A_n$ together with a line segment
connecting it to a point between $u_k$ and $l_k$. 
If $[x]$ is a diagram representing an element $x\in\SM(\A_n)$
then $\p nk(x)$ is represented by a diagram in $\D nk$ 
consisting of $\p nk'[x]$ and $k$ chords in
$\D nk\backslash L$ parallel to $L$ such that if a segment 
of $\p nk'[x]$ intersects $L$ then it 
intersects the $k$ parallel chords in $k$ over-crossings above 
the lollipop stick and $k$ under-crossings beneath the stick. See Fig \ref{f7} for a value of $\p 32$. 
We also need a linear map
$$
\g nk: \SM(\D nk)\to\SM(\D nk)
$$ 
defined by inserting 2 copies of the Jones-Wenzl idempotents $f_k$ close to $u_1,\ldots u_k$ and
$l_1,\ldots, l_k$. See Fig \ref{f7} for a value of $\g 32$. The third map
$$
\b nk : \SM(\D nk)\to \SM(\D n0)
$$
is defined by connecting $u_i$ to $l_i$ outside $\D nk$ and then pushing these arcs into $\D nk$. 
See Fig \ref{f7} for a value of $\b 32$.
The fourth map
$$
\zeta_n : \SM(\D n0) \times \SM(\A_n) \to \SM(\D 00)
$$
is defined by gluing the two entries according to the marks.
It is clear that for all $x,y\in\A_n$ we have
\begin{equation}\label{eq3}
\ph k(\la x,y\ra) = \zeta_n(\b nk\circ\g nk\circ\p nk(x),y).
\end{equation}
Therefore, if one can show that 
$\g n{k-1}\circ\p n{k-1}(b_{i_j})$, $1\le j\le s$ 
for some integer $s$, are linearly dependent then
so are the corresponding rows in $F_{n,k}$.
This observation is used to prove Lemma \ref{lemma4}.

\begin{proof}[Proof of Lemma \ref{lemma4}]
By the above argument it is enough to
show that $\g n{k-1}\circ\p n{k-1}(\bb_n)$ is contained in a subspace of dimension 
$\binom{2n}n-\binom{2n}{n-k}$ in $\SM(\D n{k-1})$. Therefore,
 it suffices to show that
\begin{equation}\label{eq4}
\dim(\im(\g n{k-1}))\le\binom{2n}n-\binom{2n}{n-k}.
\end{equation}
Let $\CL(\D n{k-1})$ be the set of non-crossing diagrams 
in $\D n{k-1}$ consisting of $n+k-1$ chords. Then $\CL(\D n{k-1})$
is a basis of $\SM(\D n{k-1})$. 
If $x\in\CL(\D n{k-1})$ contains a chord connecting two $u_i$'s or two $l_i$'s
then $\g n{k-1}(x)=0$ by a well known property of 
the Jones-Wenzl idempotent, cf. Lemma 13.2 in \cite{lickorish2}.
Hence this lemma follows from the inequality
\begin{equation*}\label{eq5}
|\tilde\CL(\D n{k-1})|\le \binom{2n}n-\binom{2n}{n-k}
\end{equation*}
where $\tilde\CL(\D n{k-1})$ is the set of diagrams in $\CL(\D n{k-1})$ with no chord connecting two $u_i$'s or two $l_i$'s.
In fact, the equality holds:
\end{proof}

\begin{lemma}\label{l4}
Assume the notation above. We have
\begin{equation}\label{eq6}
|\tilde\CL(\D nk)| = \binom{2n}n-\binom{2n}{n-k-1}.
\end{equation}
\end{lemma}

\begin{proof} Lemma 4 is a standard combinatorial fact\footnote{It can 
be derived from the reflection principle by Desir\'e Andr\'e, 1887; 
also compare \cite{Jo,SW,We,DiF}.} but we give 
its proof for completeness.
Recall that $\A_n$ denotes an annulus with $2n$ marks, 
labeled $a_1,\ldots, a_{2n}$, on the outer boundary component.
Fix a point $x_0$ between $a_{2n}$ and $a_1$ on the marked boundary 
circle such that the arc containing $x_0$ has
no other $a_i$'s.
Let $S$ be a line segment connecting $x_0$ to the other boundary 
component of $\A_n$, see \ref{fig1}..
Let $\CL_{\le k}(\A_n)\subset \bb_n$
be the set of non-crossing diagrams in $\A_n$ consisting 
of $n$ chords which intersect $S$ at most
$k$ times.
There is a 1-1 correspondence between $\CL_{\le k}(\A_n)$ and 
$\tilde\CL(\D nk)$. (Suppose $x\in\CL_{\le k}(\A_n)$
intersects $S$ at $k'$ times. Draw $k-k'$ circles close 
and parallel to the unmarked boundary component of $\A_n$. 
Cut along $S$ and we obtain a diagram in $\tilde\CL(\D nk)$.)

Hence it is enough to show that the set 
$\CL_{\ge j}(\A_n):=\bb_n\backslash\CL_{\le j-1}(\A_n)$ 
has $\binom{2n}{n-j}$ elements.  We construct a bijection 
between $\CL_{\ge j}(\A_n)$ and the choices of $n-j$ marks
among the $2n$ marks on $\A_n$.

(i) Assume $n-j$ marks, $a_{i_1},...,a_{i_{n-j}}$, are chosen. We construct 
$n$ chords as follows: If, for some $s$, the point $a_{i_s+1}$ 
is not chosen, we draw an {\em oriented} chord from $a_{i_s}$ to $a_{i_s+1}$,
 not cutting $S$. 
If $a_{2n}$ is chosen but $a_1$ is not, then draw an oriented chord 
from $a_{2n}$ to $a_1$ 
cutting $S$. At this stage at least one chord is drawn. Delete this chord 
together with its endpoint marks and 
repeat the process again until all chosen marks are used 
(they are the beginning marks of the constructed chords). We are left with $2j$ marks. 
Choose the mark with the largest index and draw a chord as before. All new $j$ 
chords cut $S$ so the constructed diagram cuts $S$ in at least $j$ points.

(ii) Conversely consider a diagram of $n$ chords cutting $S$ at least $j$ times. Orient these chords counter-clockwise.
Among the $n$ chords there are $s\le n-j$ 
of them not cutting $S$. Add to these $s$ chords $n-j-s$ more chords 
which are as close to the outside circle of $\A_n$ as possible. 
The beginning marks of these $n-j$ chords
are the marks corresponding to our diagram.

This ends the construction of the bijection. Hence we have\\ 
 $|\CL_{\ge j}(\A_n)| = \binom{2n}{n-j}$.
\end{proof}

\vspace{5ex}
\begin{minipage}[b]{0.45\linewidth}

Department of Mathematics\newline\indent
Winston-Salem State University\newline\indent
Winston Salem, NC 27110, USA\newline\indent
chenqi@wssu.edu

\end{minipage}
\hspace{3ex}
\begin{minipage}[b]{0.45\linewidth}

Department of Mathematics\newline\indent
The George Washington University\newline\indent
Washington, DC 20052, USA\newline\indent
przytyck@gwu.edu

\end{minipage}


\begin{thebibliography}{1}

\bibitem{DiF}
P.~Di Francesco, Meander determinants, {\it Comm. Math. Phys.}, 191,
1998, 543-583.

\bibitem{FK}
Igor~B. Frenkel and Mikhail~G. Khovanov.
\newblock Canonical bases in tensor products and graphical calculus for {$U\sb
  q({\mathfrak{s}}{\mathfrak{l}}\sb 2)$}.
\newblock {\em Duke Math. J.}, 87(3):409--480, 1997.

\bibitem{hoste-przytycki}
Jim Hoste and J{\'o}zef~H. Przytycki.
\newblock A survey of skein modules of {$3$}-manifolds.
\newblock In {\em Knots 90 (Osaka, 1990)}, pages 363--379. de Gruyter, Berlin,
  1992.

\bibitem{Jo}
V.~F.~R.~Jones, Index for subfactors, {\it Invent. Math.}, 72, 1983, 1-25.

\bibitem{lickorish2}
W.~B.~Raymond Lickorish.
\newblock {\em An introduction to knot theory}, volume 175 of {\em Graduate
  Texts in Mathematics}.
\newblock Springer-Verlag, New York, 1997.

\bibitem{przytycki}
J{\'o}zef~H. Przytycki.
\newblock Skein modules of {$3$}-manifolds.
\newblock {\em Bull. Polish Acad. Sci. Math.}, 39(1-2):91--100, 1991;\  
e-print: \ {\tt http://arxiv.org/abs/math/0611797}

\bibitem{schmidt}
Frank Schmidt.
\newblock Problems related to type-{$A$} and type-{$B$} matrices of chromatic
  joins.
\newblock {\em Adv. in Appl. Math.}, 32(1-2):380--390, 2004.
\newblock Special issue on the Tutte polynomial.

\bibitem{Sim} Rodica Simion,
Noncrossing partitions, {\it Discrete Math.}, 217, 2000, 367-409.

\bibitem{SW}
D.~Stanton, D.~White, in {\it Constructive Combinatorics}, Springer-Verlag,
New York-Heidelberg-Berlin, 1986.

\bibitem{We} B.~W.~Westbury,
\newblock The representation theory of the {T}emperley-{L}ieb algebras.
\newblock {\em Math. Z.}, 219(4):539--565, 1995.

\end{thebibliography}
\end{document}